\documentclass{amsart} 
\usepackage{amsmath,amscd} 
\usepackage{latexsym,amsmath,amssymb,mathrsfs,setspace,enumerate}
\usepackage[T1]{fontenc}
\usepackage[utf8]{inputenc}
\usepackage{lmodern} 
\usepackage{amssymb}
 \usepackage{tikz-cd}
\usepackage{setspace}
\usepackage{graphics}
\usepackage{graphicx}
\usepackage[utf8]{inputenc}
\usepackage[english]{babel}
\usepackage{fancyhdr}

\theoremstyle{plain}

\newtheorem{lemma}{Lemma} 
\newtheorem{prop}{Proposition} 
\newtheorem{thm}{Theorem}
\newtheorem{exam}{Example}

\newtheorem{defn}{Definition}

\newtheorem{rmk}[thm]{Remark}

\tikzcdset{%
	triple line/.code={\tikzset{%
			double distance = 3pt, 
			double=\pgfkeysvalueof{/tikz/commutative diagrams/background color}}},
	quadruple line/.code={\tikzset{%
			double distance = 5.3pt, 
			double=\pgfkeysvalueof{/tikz/commutative diagrams/background color}}},
	Rrightarrow/.code={\tikzcdset{triple line}\pgfsetarrows{tikzcd implies cap-tikzcd implies}},
	RRightarrow/.code={\tikzcdset{quadruple line}\pgfsetarrows{tikzcd implies cap-tikzcd implies}}
}
\newcommand*{\tarrow}[2][]{\arrow[Rrightarrow, #1]{#2}\arrow[dash, shorten >= 0.5pt, #1]{#2}}

\title{Category of Chain Bundles}
\author {P. G. Romeo and Riya Jose}
\address{Dept. of Mathematics, Cochin University of Science and Technology, Kochi, Kerala, INDIA.}
\email{$romeo_-parackal@yahoo.com,\, riyajosemarangattu@gmail.com $}
\thanks{} 
\date{}
\begin{document}

\begin{abstract}
		Let $\mathcal{C}$ be a category with zero. We describe the category $\mathfrak{CB}_\mathcal{C}$ of chain bundles in 
		$\mathcal{C}$ such that the objects  of  $\mathfrak{CB}_\mathcal{C}$ are 
	\[\cdots M_{3} \stackrel{Hom(M_3,M_2)}{\Rrightarrow}M_2 \stackrel{Hom(M_2,M_1)}{\Rrightarrow} M_1 \stackrel{Hom(M_1,M_0)}{\Rrightarrow} M_0 = \textbf{0} \]
	where each $M_i \in \nu \mathcal{C}$ and $\stackrel{Hom(M_{i+1},M_i)}{\Rrightarrow}$ denotes the homset 
	$Hom(M_{i+1},M_i)$. Further the homsets also includes homsets of the form $Hom(M_{i},M_i)$ and all possible composite of morphisms in $\mathcal{C}$ . The morphisms in this category are appropriate maps  between objects of 	$\mathfrak{CB}_\mathcal{C}$ and are called chain bundle maps.   Some categorical properties of 
	$\mathfrak{CB}_\mathcal{C}$ are also discussed in this paper.
\end{abstract}
\maketitle
\section{Introduction} 
Category theory has an important position in present day mathematics. It is a powerful tool which is still developing and it helps us to see the universal components of a family of structures of given kinds and how structures of different kinds are interrelated. Categories were introduced by Samuel Eilenberg and Saunders Mac Lane in the year 1945. One can formalize different branches of mathematics into categories. There are various instances where we come accross various kinds of chains in different mathematical structures where the oblects are the structures and morphisms  are structure preserving maps. In this paper we  introduce a category which we call the category of chain bundles which turns out to be interesting in several natural situations like the study of homology chains, the chains of ideals in rings and the like.

In the following we recall the basic notions and results in category theory needed in the sequel and then we introduce category of chains bundles with several examples. Further some interesting properties of this category is also discussed.
\section{Prelimanires}
We briefly recall all basic notions related to category theory, for a detailed discussion see cf.\cite{Mac}.  Subobject relation and factorization property of categories are recalled from cf.\cite{kss}.
\begin{defn} A category $\mathcal{C}$ consists of the following data:\\
	1. A class $ \nu\mathcal{C} $ called the class of vertices or objects.\\
	2. A class $\mathcal{C}$ of disjoint sets $\mathcal{C}(a,b)$, one for each pair $(a,b) \in \nu\mathcal{C} \times  \nu\mathcal{C}  $, an element $f \in \mathcal{C}(a,b)$ is called a morphism (arrow) from $a$ to $b$, written $ f: a \rightarrow b$. \\
	3. For $a,b,c \in  \nu\mathcal{C}$, a map  
	$$ \circ : \mathcal{C}(a,b) \times \mathcal{C}(b,c) \rightarrow \mathcal{C}(a,c)$$
	$$ (f,g) \rightarrow f\circ g $$\\
	is called the composition of morphisms in $\mathcal{C}$.\\
	4. For each $ a \in \nu\mathcal{C}$, a unique $1_a \in \mathcal{C}(a,a)$ is called the identity morphism on 
	$a$. These must satisfy the following axioms:\\
	
	$Cat.\, 1$ : The compostion is associative: for $f \in \mathcal{C}(a,b)$, $g \in \mathcal{C}(b,c)$ and $h \in \mathcal{C}(c,d)$ we have $ f\circ(g\circ h) = (f\circ g)\circ h.$\\
	
	$Cat.\, 2$ : For each $ a \in \nu\mathcal{C}$,  $f \in \mathcal{C}(a,b)$, $g \in \mathcal{C}(c,a), \,\,
	1_a \circ f =f $ and $ g \circ 1_a = g$\\
\end{defn}

One can identify $\nu\mathcal{C}$ as a subclass of $\mathcal{C}$, with this identification, it is possible to define categories in terms of morphisms(arrows) alone. The category $\mathcal{C}$ is said to be small if the class $\mathcal{C}$ is a set.\\ 

For any category $\mathcal{C}$ an opposite category denoted as  $\mathcal{C}^{op}$ is defined as follows:
$$ \nu \mathcal{C}^{op} = \nu \mathcal{C}, \mathcal{C}^{op}(a,b) = \mathcal{C}(b,a) \  \forall a,b \in \nu \mathcal{C}$$
and the composition $\ast$ in $\mathcal{C}^{op}$ is given by 
$$ g \ast h = h \circ g \  \forall g,h \in \mathcal{C}^{op} =\mathcal{C}$$ for which $h\circ g$ is defined.

\begin{exam}
\begin{enumerate}
		\item $\textbf{Set}$ : the category in which vertices are sets and morphisms are set maps.
		\item $\textbf{Grp}$ : the category in which vertices are groups and morphisms are homomorphisms.
		\item \textbf{Mod$_R$} : vertices are right $R$ - modules and morphisms are $R$ homomorphisms.
		\item $\textbf{Cat}$ : vertices are (small) categories and morphisms are (small) functors.
			\end{enumerate}
\end{exam}
ie., each type of algebraic systems yields a corresponding category which are called concrete categories.
\begin{defn}
	A covariant functor $F :\mathcal{C} \rightarrow \mathcal{D} $  from a category $\mathcal{C}$ to a category $\mathcal{D}$ consists of a \textit{vertex map} $ \nu F : \nu \mathcal{C} \rightarrow \nu \mathcal{D}$ which assigns to each $ a \in \nu \mathcal{C}$, a vertex $\nu F(a) \in \nu \mathcal{D}$ and a \textit{morphism map} $F : \mathcal{C} \rightarrow \mathcal{D}$ which assigns to each morphism $f: a \rightarrow b $ in $ \mathcal{C}$, a morphism $F(f) : \nu F(a) \rightarrow \nu F(b) \in \mathcal{D}$ such that \\
	$(Fn 1)$ $F(1_a) = 1_{\nu F(a)} \  \forall a \in \nu \mathcal{C}$\\
	$(Fn 2)$ $F(f)F(g) = F(fg) $ for all morphisms $f, g \in \mathcal{C}$ for which the composite $fg$ exists.
\end{defn}
	
	$F$ is a contravariant functor if $\nu F$ is as above and the morphism map assigns to each $f: a \rightarrow b $ in $ \mathcal{C}$, a morphism $F(f) : \nu F(b) \rightarrow \nu F(a) \in \mathcal{D}$ such that they satisfy axiom $(Fn 1)$ and $(Fn^* 2)$ ie.,  $F(g)F(f) = F(fg) $ for all morphisms $f, g \in \mathcal{C}$ for which the composite $fg$ exists.

\begin{exam}
	A category $\mathcal{D}$ is a subcategory of a category $\mathcal{C}$ if the class $\mathcal{D}$ is a subclass of $\mathcal{C}$ and the composition in $\mathcal{D}$ is the restriction of the composition in $\mathcal{C}$ to $\mathcal{D}$. In this case, the inclusion $\mathcal{D} \subseteq \mathcal{C}$ preserves composition and identities and so represents a functor of $\mathcal{D}$ to $\mathcal{C}$ which is called the inclusion functor of $\mathcal{D}$ into $\mathcal{C}$.
\end{exam}
\begin{defn} 
	A morphism $f$ in a category $\mathcal{C}$ is a monomorphism if $ gf = hf \Rightarrow g = h \ \forall g,h \in \mathcal{C} $, 	i.e., $f$ is a monomorphism if it is right cancellable. A morphism $f$ in  $\mathcal{C}$ is called a split monomorphism if it has a right inverse.
\end{defn}	
	
	Every morphism in a concrete category whose underlying function is an injection is a monomorphism. 	A morphism 
	$f$ in a category $\mathcal{C}$ is an epimorphism if $ fg = fh \Rightarrow g = h \ \forall g,h 
	\in \mathcal{C}$, ie., $f$ is en epimorphism if it is left cancellable. A morphism $f$ in  $\mathcal{C}$ is called a split epimorphism if it has a left inverse. Every morphism in a concrete category whose underlying function is an surjection is an epimorphism.

\begin{defn}
	A Preorder $\mathcal{P}$ is a category such that for any $p, p' \in  \mathcal{P}$, $\mathcal{P}(p,p')$ contains atmost one morphism. In this case, the relation $\subseteq$ on the class $\nu \mathcal{P}$ defined by 
	\begin{equation}
	p \subseteq p' \ \Leftrightarrow \ \mathcal{P}(p,p') \neq \emptyset
	\end{equation}
	is a quasiorder. When $\mathcal{P}$ is a preorder, $\nu \mathcal{P}$ will stand for the quasiordered class $(\nu \mathcal{P}, \subseteq)$. Conversely given a quasiorder $\leq$ on the class $X$, the subset $$ \mathcal{P} = \{(x,y) =\in X \times X : x \leq y\}$$
	of $ X \times X $ is a preorder such that the quasiordered class $\nu \mathcal{P}$ defined above is order isomorphic with $(X,\leq)$
\end{defn}	

	If the relation $\subseteq$ on $\mathcal{P}$  is antisymmetric then we shall say that $\mathcal{P}$ is a 
	strict preorder.

\begin{defn} 
	Let $ \mathcal{C} $ be a category. A choice of subobjects in $ \mathcal{C} $ is a subcategory $ \mathcal{P} \subseteq \mathcal{C} $ satisfying the following:\\
	\begin{enumerate}
		\item[a.] $\mathcal{P}$ is a strict preorder with $\nu \mathcal{P} = \nu \mathcal{C}$
		\item[b.] Every $f \in \mathcal{P} $ is a monomorphism in  $\mathcal{C}$
		\item[c.] If $f, g \in \mathcal{P}$ and if $f = hg$ for some $h \in \mathcal{C}$ then $h \in \mathcal{P}$.
	\end{enumerate}
	when $\mathcal{P}$ satisfies these conditions, the pair $ (\mathcal{C},\mathcal{P})$ is called a caegory with subobjects.
\end{defn}

\begin{rmk}
	When $\mathcal{C} $ has  subobjects, unless explicitly stated otherwise, $ \nu\mathcal{C} $ will denote the choice of subobjects in $ \mathcal{C} $. The partial order defined by equation $(1)$ is called the preorder of inclusions or subobject relation in $ \mathcal{C} $ and is denoted by $ \subseteq $. If $c,d \in \nu\mathcal{C}$ and $ c\subseteq d $ the unique morphism from $c$ to $d$ is the inclusion $j_c^d: c \rightarrow d $.
\end{rmk}
\begin{exam}
	In categories $Set, Grp, Vct_K, Mod_R $ the relation on objects induced  by the usual set inclusion is a subobject relation.
\end{exam}

\begin{defn}
	A morphism $f$ in a category $\mathcal{C}$ with subobjects is said to have factorization if $f$ can be expressed as  $ f = pm$ where $p$ is an epimorphism and $m$ is an embedding.
\end{defn}

	Factorization of a morphism need not be unique. Every morphism $f$ with factorization has atleast one factorization of the form $f = qj$  where $q$ is an epimorphism and $j$ is an inclusion, such factorizations are called canonical factorization. $\mathcal{C}$ is a category with factorization if $\mathcal{C}$ has subobjects and if every morphism in $\mathcal{C}$ has factorization. The category has unique factorization property if every morphism in $\mathcal{C}$ has unique canonical factorization.

\begin{exam}
	If $f : X \rightarrow Y$ is a mapping of sets and $f(X) = Im f$ then $f(X) \subseteq  Y$ and we can write $f = f^0j_{f(X)}^Y$. Here $f^0$  denote the mapping of $X$ onto $f(X)$ determined by $f$. Since surjective mappings are epimorphisms in $Set$, this gives a canonical factorization of $f$ in $Set$ which is clearly unique. Thus $Set$ is a category with unique factorization.
\end{exam}	
	Since surjective continous mappings are epimorphisms in $Top$, it follows as in the above example that this category has factorization property. However, if $Y$  is dense in $X$, $ h =j_Y^X$ is an epimorphism in $Top$ and $ h = 1_Yj_Y^X = j_Y^X1_X$. Then both $1_Yj_Y^X $ and $ j_Y^X1_X$ are canonical factorizations of $h$ in $Top$. Thus $Top$ doesnot have unique factorization property.

\begin{prop}
	Let $\mathcal{C}$ be category with factorization. Suppose that the morphism $f \in \mathcal{C}$ has the following property:\\
	$(Im)$ $f$ has a canonical factorization $f = xj$ such that for any canonical factorization $f = yj'$ of $f$, there is an inclusion $j''$ with $y = xj''$. Then the factorization $f = xj$ is unique.
\end{prop}

\begin{rmk}
	A morphism $f$ in a category with factorization is said to have \textit{image} if $f$ satisfies the condition $(Im)$ of the Proposition above. In this case the unique canonical factorization $f = xj$ with the property $(Im)$ is denoted by $ f = f^0j_f$ where $f^0$ is called the \textit{epimorphic component} of $f$ and $j_f$ is called the \textit{inclusion} of $f$. The unique vertex $Im f = cod f^0 = dom j_f$ is called the \textit{image} of $f$.
\end{rmk}

Since categories $Set, Grp, $etc., has unique factorization, morphisms in these categories have images. Though the category $Top$ doesnot have unique factorization, it can be seen that every morphism in $Top$ also has image.

\begin{defn}
	A groupoid $\mathcal{G}$ is a small category in which every morphism is an isomorphism. A groupoid  $\mathcal{G}$ is said to be connected if for all $a \in \nu \mathcal{G}, H_a = \mathcal{G}(a,a) \neq \phi$.
\end{defn}

\begin{exam}
	Every group $G$ is a groupoid with exactly one vertex.
\end{exam}

\section{Category of Chain Bundles}
Let $\mathcal{C}$ be a category  with zero. Now we describe the category of chain bundles  and cochain bundles and it is shown that they are in fact categories with subobjects and admits subcategories with factorization property.
Some examples of such categories are also provided. 
\begin{defn}
	Let $\mathcal{C}$ be category with zero. A chain bundle in $\mathcal{C}$ is of the form 
	\[\cdots M_{3} \stackrel{Hom(M_3,M_2)}{\Rrightarrow}M_2 \stackrel{Hom(M_2,M_1)}{\Rrightarrow} M_1 \stackrel{Hom(M_1,M_0)}{\Rrightarrow} M_0 = \textbf{0} \]
	where vertices $M_i \in \nu \mathcal{C}$ for all $i$ and morphisms are homsets $Hom(M_{i+1},M_i)$ denoted by 
	$\stackrel{Hom(M_{i+1},M_i)}{\Rrightarrow}$. The chain bundle consists of all homstes of the form 
	$Hom(M_{i},M_i)$ and all possible composite of morphisms. 
\end{defn}	

	Note that a chain bundle in category $\mathcal{C}$ is a subcategory of $\mathcal{C}$. A cochain bundle in a category 
	$\mathcal{C}$ is a subcategory of the form 
	\[ \textbf{0} = M_0\stackrel{Hom(M_0,M_1)}{\Rrightarrow}  M_{1} \stackrel{Hom(M_1,M_2)}{\Rrightarrow}M_2 \stackrel{Hom(M_2,M_3)}{\Rrightarrow} M_3  \cdots \]

\begin{defn}
	A chain bundle map between two chain bundles is a functor $F$ between the two, such that resulting diagram commutes and the vertex map is a sequence of morphisms in $\mathcal{C}$, $\nu F = \{f_i: M_i \rightarrow N_i\} $ and morphism map is a map on homset of $\mathcal{C}$
	\begin{center}
		\begin{tikzcd}
		\cdots \tarrow[" "]{r} & M_3 \tarrow[" "]{r} \arrow{d}[swap]{f_3}
		& M_2 \tarrow[" "]{r} \arrow{d}{f_2} & M_1 \arrow{d}{f_1} \tarrow[" "]{r} & \textbf{0} \arrow{d}{f_0}\\
		\cdots \tarrow[" "]{r}& N_3 \tarrow [" "]{r}  & N_2 \tarrow[" "]{r} &  N_1 \tarrow[" "]{r} & \textbf{0}  \\
		\end{tikzcd}\\
	\end{center}
\end{defn}	

	Similarly we define cochain bundle map. A category of chain bundles is the one whose objects are chain bundles in $\mathcal{C}$ and morphism between two chain bundles are chain bundle maps. We denote this category as 
	$\mathfrak{CB}_{\mathcal{C}}$

\begin{exam}
	Let $\mathcal{C} $ be the category of $\mathbb{Z}$- modules and $ \mathfrak{CB}$ be its category of chain bundles. Consider two chains bundles $$ 3\mathbb{Z} \stackrel{\frac{2}{3} a} {\Rrightarrow} 2\mathbb{Z} \stackrel{\frac{5}{2} b}{\Rrightarrow} 5\mathbb{Z}\stackrel{0}{\Rrightarrow} \textbf{0} \quad \text{and} \quad  6\mathbb{Z}\stackrel{\frac{2}{3} a'}{\Rrightarrow} 4\mathbb{Z} \stackrel{\frac{1}{2} b'}{  \Rrightarrow} 2\mathbb{Z} \stackrel{0}{\Rrightarrow} \textbf{0}$$ 
	A morphism $F$ between these two chain bundles is as follows:

\begin{center}
	\begin{tikzcd}
3 \mathbb{Z} \tarrow["\frac{2}{3} a"]{r} \arrow{d}[swap]{2}
& 2 \mathbb{Z} \tarrow["\frac{5}{2} b"]{r} \arrow{d}{2} & 5 \mathbb{Z} \arrow{d}{\frac{1}{5}} \tarrow["0"]{r} & \textbf{0} \arrow{d}{0}\\
6 \mathbb{Z} \tarrow ["\frac{2}{3} a'"]{r}  & 4 \mathbb{Z} \tarrow["\frac{1}{4} b'"]{r} &  \mathbb{Z} \tarrow["0"]{r} & \textbf{0}  \\
\end{tikzcd}
\end{center}
\flushleft{ $\nu F(3\mathbb{Z}) = 6\mathbb{Z} , \nu F(2\mathbb{Z}) = 4\mathbb{Z} , \nu F(5\mathbb{Z}) = 2\mathbb{Z},\nu F(\textbf{0}) = \textbf{0} $} and $F(\frac{2}{3} a) = \frac{2}{3} a',\, F(\frac{5}{2} b) = \frac{1}{4}b', 
\,F(0) = 0.  $\\
Another morphism $G$ between these to chain bundles is as follows:

\begin{center}
	\begin{tikzcd}
	3 \mathbb{Z} \tarrow["\frac{2}{3} a"]{r} \arrow{d}[swap]{4}
	& 2 \mathbb{Z} \tarrow["\frac{5}{2} b"]{r} \arrow{d}{4} & 5 \mathbb{Z} \arrow{d}{\frac{2}{5}} \tarrow["0"]{r} & \textbf{0} \arrow{d}{0}\\
	6 \mathbb{Z} \tarrow ["\frac{2}{3} a'"]{r}  & 4 \mathbb{Z} \tarrow["\frac{1}{4} b'"]{r} &  \mathbb{Z} \tarrow["0"]{r} & \textbf{0}  \\
	\end{tikzcd}
\end{center}
\flushleft{ $\nu G(3\mathbb{Z}) = 6\mathbb{Z} , \nu G(2\mathbb{Z}) = 4\mathbb{Z} , \nu G(5\mathbb{Z}) = 2\mathbb{Z},\nu G(\textbf{0}) = \textbf{0} $} and $G(\frac{2}{3} a) = \frac{2}{3} a',\, G(\frac{5}{2} b) = \frac{1}{4}b',\, G(0) = 0.  $\\
\end{exam}

\begin{rmk}
	The two morphisms in the above example are same when considered as functors but are different when considered as chain bundle maps.
\end{rmk}
\begin{exam}
	Let $\mathcal{C}$ be the category of subgroups of symmetric group $S_3$, and $ \mathfrak{CB} $ the category of chain bundles in $\mathcal{C}$. Consider the two chain bundles 
	$$ S_3 \stackrel{0} {\Rrightarrow} A_3 \stackrel{0}{\Rrightarrow}  \textbf{0} \quad \text{and} \quad  K_1 \stackrel{Hom(K_1,S_3)}{\Rrightarrow} S_3 \stackrel{0}{\Rrightarrow} \textbf{0}$$ 
	$K_1=\{ \textbf{0}, (2,3)\}$, then one of the morphism $F$ between these to chain bundles is as follows:
	
\begin{center}
		\begin{tikzcd}
		S_3 \tarrow["0"]{r} \arrow{d}[swap]{0}
		& A_3 \tarrow["0"]{r} \arrow{d}{p} &  \textbf{0} \arrow{d}{0}\\
		K_1 \tarrow ["{Hom(K_1,S_3)} "]{r}  & S_3  \tarrow["0"]{r} & \textbf{0}  \\
		\end{tikzcd}
	\end{center}
	
\flushleft{ $\nu F(S_3) = K_1 , \nu F(A_3) =S_3 , \nu F(\textbf{0}) =  \textbf{0} $}. $p$ can be any morphism from $A_3$ to $S_3$. The map $ 0 : S_3 \rightarrow A_3 $ can be mapped to any  $ g: K_1 \rightarrow S_3$ .
	\end{exam}
	
\begin{defn}
	Let $\mathcal{C}$ be a category with subobjects. Let $M_1, M_2, M_{1}',M_2' \in \nu \mathcal{C}$. A morphism $f : M_1' \rightarrow M_2' $ is said to be corestriction of a morphism $g: M_{1} \rightarrow M_2 $ to inclusion $i: M_2' \rightarrow M_2 $ if $M_{1}' \subset M_1 , M_2' \subset M_2$ and $(g|_{M_{1}'})^0 = f$. 
\end{defn}

\begin{defn}
	A chain bundle $c'$:\[\cdots M_{3}' \stackrel{Hom(M_3',M_2')}{\Rrightarrow}M_2' \stackrel{Hom(M_2',M_1')}{\Rrightarrow} M_1' \stackrel{Hom(M_1',M_0')}{\Rrightarrow} M_0' = \textbf{0} \]
	is a subchain bundle of the chain bundle $c$: \[\cdots M_{3} \stackrel{Hom(M_3,M_2)}{\Rrightarrow}M_2 \stackrel{Hom(M_2,M_1)}{\Rrightarrow} M_1 \stackrel{Hom(M_1,M_0)}{\Rrightarrow} M_0 = \textbf{0} \]
	if $M_i'$ is a subobject of $M_i$ and for each $f' \in Hom(M_{i+1}',M_i') $ there is a $f \in Hom(M_{i+1},M_i)$ such that $f'$ is corestriction of $f$.
\end{defn}	
If $c'$ is a subchain bundle of $c$ then there exists a chain bundle map from $c'$ to $c$ whose each vertex map is an inclusion and morphism map on each homset is such that each morphism in  $Hom(M_{i+1}',M_i')$ is a corestriction of some  morphism in $Hom(M_{i+1},M_i)$.

\begin{center}
	\begin{tikzcd}
	\cdots \tarrow[" "]{r} & M_3' \tarrow[" "]{r} \arrow{d}[swap]{i}
	& M_2' \tarrow[" "]{r} \arrow{d}{i} & M_1' \arrow{d}{i} \tarrow[" "]{r} & \textbf{0} \arrow{d}{i}\\
	\cdots \tarrow[" "]{r}& M_3 \tarrow [" "]{r}  & M_2 \tarrow[" "]{r} &  M_1 \tarrow[" "]{r} & \textbf{0}  \\
	\end{tikzcd}\\
\end{center}

\begin{lemma}
	The realtion \textquotedblleft is a subchain bundle" is a strict preorder on vertex class of category of chain bundles.
\end{lemma}

\begin{lemma}
	If $\mathcal{C}$ is a category with subobjects, then category of chain bundles in $\mathcal{C}$, $\mathfrak{CB}_\mathcal{C}$ forms a category with subobjects where being subchain bundle is the subobject relation.
\end{lemma}

\begin{exam}
	Consider the chain bundle $c: 3\mathbb{Z} \Rrightarrow 2\mathbb{Z} \Rrightarrow 8\mathbb{Z} \Rrightarrow \textbf{0}$ in category of chain bunbles in category of submodules of $\mathbb{Z}$. Trivally $c$ is a subobject of itself.\\
	\begin{center}
		\begin{tikzcd}
		3 \mathbb{Z} \tarrow["\frac{2}{3} a"]{r} \arrow{d}[swap]{i}
		& 2 \mathbb{Z} \tarrow["4b"]{r} \arrow{d}{i} & 8 \mathbb{Z} \arrow{d}{i} \tarrow["0"]{r} & \textbf{0} \arrow{d}{i}\\
		3 \mathbb{Z} \tarrow["\frac{2}{3} a"]{r} 
		& 2 \mathbb{Z} \tarrow["4b"]{r} & 8 \mathbb{Z} \tarrow["0"]{r} & \textbf{0} \\
		\end{tikzcd}
	\end{center}
	
	Further it is also seen that $c':6\mathbb{Z} \Rrightarrow 4\mathbb{Z} \Rrightarrow 16\mathbb{Z} \Rrightarrow \textbf{0}$ is a subobject of $c$.
	
	\begin{center}
		\begin{tikzcd}
		6 \mathbb{Z} \tarrow["\frac{2}{3} a'"]{r} \arrow{d}[swap]{i}
		& 4 \mathbb{Z} \tarrow["4b'"]{r} \arrow{d}{i} & 16 \mathbb{Z} \arrow{d}{i} \tarrow["0"]{r} & \textbf{0} \arrow{d}{i}\\
		3 \mathbb{Z} \tarrow["\frac{2}{3} a"]{r} 
		& 2 \mathbb{Z} \tarrow["4b"]{r} & 8 \mathbb{Z} \tarrow["0"]{r} & \textbf{0} \\
		\end{tikzcd}\\
	\end{center}
where morphism map is identity on each homsets with obvious inclusion from $c'$ to $c$. Where as $c'':9\mathbb{Z} \Rrightarrow 4\mathbb{Z} \Rrightarrow 16\mathbb{Z} \Rrightarrow \textbf{0}$ is a not a subobject of $c$. Though each vertex of $c''$ is a subobject of $c$, the morphisms in  $Hom(9 \mathbb{Z}, 4 \mathbb{Z})$ is not restriction of that in  $Hom(3 \mathbb{Z}, 2 \mathbb{Z})$. 
	
	\begin{center}
		\begin{tikzcd}
		9 \mathbb{Z} \tarrow["\frac{4}{9} a'"]{r} \arrow{d}[swap]{i}
		& 4 \mathbb{Z} \tarrow["4b'"]{r} \arrow{d}{i} & 16 \mathbb{Z} \arrow{d}{i} \tarrow["0"]{r} & \textbf{0} \arrow{d}{i}\\
		3 \mathbb{Z} \tarrow["\frac{2}{3} a"]{r} 
		& 2 \mathbb{Z} \tarrow["4b"]{r} & 8 \mathbb{Z} \tarrow["0"]{r} & \textbf{0} \\
		\end{tikzcd}\\
	\end{center}
\end{exam}

Let $\mathcal{C}$ be category with factorization. Then it is possible to obtain a subcategory of  $\mathfrak{CB}_\mathcal{C}$ denoted by $\mathfrak{CB}_\mathcal{C}'$  whose vertex class is same as that of $\mathfrak{CB}_\mathcal{C}$ but morphisms are those chain bundle maps  in $\mathfrak{CB}_\mathcal{C}$ whose morphism map is full. 
Let $F$ be a chain bundle map in $\mathfrak{CB}_\mathcal{C}'$. Then $F$ can be factorized as follows:

$$  \cdots M_{3} {\Rrightarrow}M_2 {\Rrightarrow} M_1 {\Rrightarrow}\textbf{0} $$
\hspace{7.5cm} 	$ \downarrow $ $F^0$ epimorphism

$$\cdots N_{3}' {\Rrightarrow}N_2' {\Rrightarrow} N_1' {\Rrightarrow} \textbf{0}$$
\hspace{7.5cm}	$ \downarrow $ inclusion

$$\cdots N_{3} {\Rrightarrow}N_2 {\Rrightarrow} N_1 {\Rrightarrow} \textbf{0}$$	
where $N_i' = F(N_i)$. Thus $\mathfrak{CB}_\mathcal{C}'$ forms a category with factorization.

\begin{exam}
	Consider the following two chain bundles and functor $F$ between them in category of chain bundles in category of submodules of $\mathbb{Z}$.\\
	
	\begin{center}	
		\begin{tikzcd}
		3\mathbb{Z} \tarrow["{\frac{2}{3} a}"]{r} \arrow{d}[swap]{4}
		& 2\mathbb{Z} \tarrow["{\frac{5}{2} b}"]{r} \arrow{d}{4} & 5\mathbb{Z} \arrow{d}[swap]{\frac{2}{5}}\tarrow["0"]{r} & \textbf{0} \arrow{d}{0}\\
		6\mathbb{Z} \tarrow["{\frac{2}{3} a'}"]{r}
		& 4\mathbb{Z} \tarrow["{\frac{1}{4} b'}"]{r}  & \mathbb{Z} \tarrow["0"]{r} & \textbf{0} \\
		\end{tikzcd}	
	\end{center}
	\flushleft{ $\nu F(3\mathbb{Z}) = 6\mathbb{Z} , \nu F(2\mathbb{Z}) =4 \mathbb{Z} , \nu F(5\mathbb{Z}) = \mathbb{Z}$} and $F(\frac{2}{3} a) = \frac{2}{3} a, F(\frac{5}{2} b) = \frac{1}{4} b,$ ie., $F$ can be factorized into 
	$G \circ H$ where $G$ is an epimorphism and $H$ is an inclusion. Factorization is as follows:	
	
	\begin{center}
		\begin{tikzcd}
		3\mathbb{Z} \tarrow["{\frac{2}{3} a}"]{r} \arrow{d}[swap]{4}
		& 2\mathbb{Z} \tarrow["{\frac{5}{2} b}"]{r} \arrow{d}{4} & 5\mathbb{Z} \arrow{d}[swap]{\frac{2}{5}}\tarrow["0"]{r} & \textbf{0} \arrow{d}{0}\\
		12\mathbb{Z} \tarrow["{\frac{2}{3} a'}"]{r} \arrow{d}[swap]{i}
		& 8\mathbb{Z} \tarrow["{\frac{1}{4} b'}"]{r} \arrow{d}{i} & 2\mathbb{Z} \arrow{d}[swap]{i}\tarrow["0"]{r} & \textbf{0} \arrow{d}{i}\\
		6\mathbb{Z} \tarrow["{\frac{2}{3} a'}"]{r}
		& 4\mathbb{Z} \tarrow["{\frac{1}{4} b'}"]{r}  & \mathbb{Z} \tarrow["0"]{r} & \textbf{0} \\
		\end{tikzcd}
	\end{center}
	\flushleft{ $\nu G(3\mathbb{Z}) = 12\mathbb{Z} , \nu G(2\mathbb{Z}) = 8\mathbb{Z} , \nu G(5\mathbb{Z}) = 2\mathbb{Z}$} and $G(\frac{2}{3} a) = \frac{2}{3} a, G(\frac{5}{2} b) = \frac{1}{4} b,$\\ 
	\flushleft{ $\nu H(12\mathbb{Z}) = 6\mathbb{Z} , \nu H(8 \mathbb{Z}) = 4 \mathbb{Z} , \nu H(2\mathbb{Z}) = \mathbb{Z}$} and 	$H(\frac{2}{3} a') = \frac{2}{3} a', H(\frac{1}{4} b') = \frac{1}{4} b'$\\ 
\end{exam}

\begin{rmk}
	Similarly we can construct category of cochain bundles $\mathfrak{CCB}$ and it can be show that it is a category with subobjects.
\end{rmk}

\subsection{Groupoids from Category of Chain Bundles}
If $\mathcal{G}$ is a groupoid, then the category $\mathfrak{CB}_\mathcal{G}$ has a subcategory $\mathfrak{CB}_\mathcal{G}'$ which form a groupoid.  $\mathfrak{CB}_\mathcal{G}'$ has same vertex class as that of $\mathfrak{CB}_\mathcal{G}$ but we choose only those morphisms whose morphism map is as follows:

\begin{center}
	\begin{tikzcd}
	\cdots \tarrow[" "]{r} & M_3 \tarrow[" "]{r} \arrow{d}[swap]{f_3}
	& M_2 \tarrow[" "]{r} \arrow{d}{f_2} & M_1 \arrow{d}{f_1} \tarrow[" "]{r} & \textbf{0} \arrow{d}{f_0}\\
	\cdots \tarrow[" "]{r}& N_3 \tarrow [" "]{r}  & N_2 \tarrow[" "]{r} &  N_1 \tarrow[" "]{r} & \textbf{0}  \\
	\end{tikzcd}\\
\end{center}
 Choose $F$ as a morphism in $\mathfrak{CB}_\mathcal{G}'$ only if its morphism map is of the form $F(f) = f_i^{-1}ff_{i-1}^{-1}$ for each $f \in Hom(M_i,M_{i-1})$.\\
 
 Products and coproducts in the category $\mathfrak{CB}$ are done termwise. If two chain bundles have different lengths then we add zeros as vertices to left of chain bundle to make length of both equal.\\
 Consider the following two chain bundles $c$ and $d$:\\
 \[c:\cdots M_{3}' {\Rrightarrow}M_2' {\Rrightarrow} M_1' {\Rrightarrow} M_0' = \textbf{0} \]
 \[d: \cdots M_{3} {\Rrightarrow}M_2 {\Rrightarrow} M_1 {\Rrightarrow} M_0 = \textbf{0} \] 
 The product of $c \times d $ is the chain bundle map 
  \[c \times d:\cdots M_{3}' \times M_{3} {\Rrightarrow}M_2'\times M_{2} {\Rrightarrow} M_1'\times M_{1}  {\Rrightarrow} M_0'\times M_{0}  = \textbf{0} \times \textbf{0}  \]
  if for any $F: l \rightarrow c$ and $G: l \rightarrow d $ where $l$ is the chain bundle 
  $$l:\cdots N_{3} {\Rrightarrow}N_2 {\Rrightarrow} N_1 {\Rrightarrow} N_0 = \textbf{0}  $$
  then there exists a chain bundle map $L : l \rightarrow c \times d $ such that for any $k \in Hom(N_i,N_j) $ and $L(k) \in Hom(M_i' \times M_i, M_j' \times M_j) $ there corresponds a $F(k) \in Hom(M_i' \times M_j')$ and $G(k) \in Hom(M_i,M_j)$
  
\section{Category of Chains}
Let $\mathcal{C}$ be a category. $\mathfrak{CB}_\mathcal{C}'$ be the subcategory of the chain bundle category 
$\mathfrak{CB}_\mathcal{C}$ whose morphism map is full and let 
\[\cdots M_{3} \stackrel{Hom(M_3,M_2)}{\Rrightarrow}M_2 \stackrel{Hom(M_2,M_1)}{\Rrightarrow} M_1 \stackrel{Hom(M_1,M_0)}{\Rrightarrow} M_0 = \textbf{0} \]  be a chain bundle in $\mathfrak{CB}_\mathcal{C}'$. Then by choosing atmost one morphism from each homset in the chain bundle, we obtain  chains in $\mathfrak{CB}_\mathcal{C}'$ . 

\begin{exam}
	Consider the chain bundle $$18 \mathbb{Z} \Rrightarrow 9 \mathbb{Z} \Rrightarrow 3\mathbb{Z}\Rrightarrow 8 \mathbb{Z} \Rrightarrow 4 \mathbb{Z} \Rrightarrow   2\mathbb{Z} \Rrightarrow  5\mathbb{Z}\Rrightarrow  \textbf{0}$$
	by choosing only inclusions from each homset we get following chains:\\
	$18 \mathbb{Z} \rightarrow 9 \mathbb{Z} \rightarrow 3\mathbb{Z}\rightarrow \textbf{0},\,\,8 \mathbb{Z} \rightarrow 4 \mathbb{Z} \rightarrow   2\mathbb{Z} \rightarrow \textbf{0}, 5\mathbb{Z}\rightarrow  \textbf{0}$
\end{exam}
Now define chain maps as follows:  
\begin{defn}
	A chain map between two chains as a sequence of morphisms in $\mathcal{C},  \{f_i: N_i \rightarrow M_i\} $ such that the resulting diagram commutes.
	\begin{center}
		\begin{tikzcd}
		\cdots \arrow[" "]{r} & N_3 \arrow[" "]{r} \arrow{d}[swap]{f_3}
		& N_2 \arrow[" "]{r} \arrow{d}{f_2} & N_1 \arrow{d}{f_1} \arrow[" "]{r} & \textbf{0} \arrow{d}{f_0}\\
		\cdots \arrow[" "]{r}& M_3 \arrow [" "]{r}  & M_2 \arrow[" "]{r} &  M_1 \arrow[" "]{r} & \textbf{0}  \\
		\end{tikzcd}\\
	\end{center}
\end{defn} 

Given a category $\mathcal{C}$, it is seen that one can construct category of chain bundles $ \mathfrak{CB}_{\mathcal{C}}$. Now we choose another category from  $ \mathfrak{CB}_{\mathcal{C}}$ by choosing a morphism from each homset to form chains. Depending on the choice we get different collection of chains from $ \mathfrak{CB}_{\mathcal{C}}$. Let $\nu \Gamma$ be one such collection of chains obtained from $ \mathfrak{CB}_{\mathcal{C}}$, 
we obtain a category  $\Gamma$ whose vertex class is $\nu \Gamma$ and morphisms between two chains is chain map defined above.

\begin{exam}
	Let $\mathcal{C}$ be an abelian category. Consider the category $ \mathfrak{CB}_{\mathcal{C}}$ of chain bundles in $\mathcal{C}$. We set the condition for choosing chain from a chain bundle in $ \mathfrak{CB}_{\mathcal{C}}$. 
	Let  $c: \cdots C_{n+1} \Rrightarrow C_{n} \Rrightarrow C_{n-1}\Rrightarrow \cdots \Rrightarrow \textbf{0} $ be a chain bundle in  $ \mathfrak{CB}_{\mathcal{C}}$, we choose one  $\partial_i  $ from each homset $Hom (C_{i+1},C_{i})$ such that $\partial_{i+1}\circ \partial_{i} = 0 $. If we can choose such a $\partial_i  $ from each homset $Hom (C_i,C_{i+1}) \quad \forall i$ we get a single chain of the form 
$$c_1: \cdots C_{n-1} \stackrel{\partial_{n-1}}{\rightarrow} C_n \stackrel{\partial_n}{\rightarrow} C_{n+1} \stackrel{\partial_{n+1}}{\rightarrow} \cdots .$$
	otherwise we get a collection of such chains from the single chain bundle. Thus we can construct a category $\Gamma$ whose vertices are chains obtained from chain bundles in $\mathfrak{CB}_{\mathcal{C}}$ as above and morphisms between two chains is defined as sequence of morphisms in $\mathcal{C}$ from vertices of first chain to that of second one which make the resulting diagram commutative. The category of chain complexes in $\mathcal{C}$ then coincides with the category  $\Gamma$.
\end{exam}
If $\mathcal{C}$ is a category with subobjects, then so is $\mathfrak{CB}_{\mathcal{C}}$ and similarly category of chains obtained from it is also a category with subobjects.


\begin{thebibliography}{99}

\bibitem{Mac} Saunders Mac Lane: Categories for the Working Mathematician, Second edition, 0-387-9803-8, Springer-Verlag New york, Berlin Heidelberg Inc., 1998.
\bibitem{kss} K.S.S. Nambooripad : Theory of Regular Semigroups, Sayahna Foundation Trivandrum, 2018.

\bibitem{berrik} A. J. Berrick and M. E. Keating: Categotries and Modules, Cambridge University Press, 2000. 

\bibitem{david}  Davide L. Ferrario, Renzo A. Piccinini: Simplicial Structures in Topology.Springer Science and Business Media, 2010.

\end{thebibliography}
\end{document}